\documentclass[11pt,a4paper]{article}
\usepackage[margin=2.5cm]{geometry}
\usepackage{amsmath,amssymb}
\usepackage{graphicx}
\usepackage{booktabs}
\usepackage{microtype}
\usepackage{hyperref}
\newcommand{\E}{\mathbb{E}}
\newcommand{\eps}{\varepsilon}

\title{Monte Carlo pricing under fast mean-reverting stochastic volatility:\\
the multi-scale limit $\eps\to0$}
\author{ Laurent Mertz (Math. dept. City University, Hong Kong, China) \footnote{fmertz@cityu.edu.hk},\and Olivier Pironneau (applied math. dept LJLL, Sorbonne University, Paris, France)\footnote{olivier.pironneau@sorbonne-universite.fr}.} 
\date{\today}

\begin{document}
\maketitle

\begin{abstract}
We compute $\E[(S_T-K)^+]$ by Monte Carlo for a scalar stochastic-volatility model with a fast
mean-reverting factor of time scale $\eps$, for $\eps$ ranging from $1$ down to $10^{-3}$.
A conditional (mixing) estimator gives finite variance, whereas the direct estimator has infinite
variance for this model. The volatility factor is simulated with its exact
Ornstein--Uhlenbeck transition. As $\eps\to0$ the price converges, at rate $O(\eps)$, to the
Black--Scholes price with the averaged volatility $\bar\sigma$, and the implied-volatility smile flattens to $\bar\sigma$.
Finally, we test a martingale control variate built on the Black--Scholes delta with volatility $\bar\sigma$.
If the martingale is driven by the true volatility $\sigma(Y_t)$, the variance is reduced by a factor that grows like $1/\eps$,
about $160$ at $\eps=10^{-3}$. If it is driven by the constant $\bar\sigma$, the variance is essentially not reduced.
Last, we calibrate the model with $\rho\neq0$ to S\&P~500 implied volatilities by full simulation of $S$.
At the same number of paths the control variate reduces the variance by a factor $1.7$--$16$ (median $5.6$)
and gives more accurate calibrated parameters; at equal CPU time it pays off only if the delta is rebalanced
on a coarser grid than the time step.
For a basket of four indices (state dimension $8$) the gain is larger (median $8$) and the relative cost smaller,
so that the control variate is about $4$--$5$ times more efficient than plain Monte Carlo at equal CPU time.
\end{abstract}

\section{Model}
Let $W^0,W^1$ be independent Brownian motions ($\rho=0$) and
\begin{equation}\label{eq:model}
dS_t = rS_t\,dt + S_t\,\sigma(Y_t)\,dW^0_t,\qquad
dY_t = \frac1\eps(m_1-Y_t)\,dt + \frac{\nu\sqrt2}{\sqrt\eps}\,dW^1_t,\qquad
\sigma(y)=e^{y+\theta},
\end{equation}
with $r=0.08$, $m_1=-1$, $\nu=1$, $\theta=0.5$, $S_0=100$, $Y_0=-1$, $T=1$.
We compute $\E[(S_T-K)^+]$ (undiscounted); the price is $e^{-rT}$ times this value.

The process $Y$ is an Ornstein--Uhlenbeck process with rate $1/\eps$ and invariant law
$\mathcal N(m_1,\nu^2)$, independent of $\eps$. When $\eps\to0$, ergodicity gives
\[
\frac1T\int_0^T\sigma^2(Y_t)\,dt \;\longrightarrow\; \bar\sigma^2
 = \langle\sigma^2\rangle = e^{2\theta}\,\E\big[e^{2Y}\big]
 = e^{2\theta+2m_1+2\nu^2},
\qquad \bar\sigma = e^{0.5}\approx 1.6487,
\]
so the limit price is the Black--Scholes price with constant volatility $\bar\sigma$.

\section{Method}
\paragraph{Direct simulation has infinite variance.}
Conditionally on $Y$, $\E[S_T^2\mid Y]=S_0^2e^{2rT}e^{Q}$ with $Q=\int_0^T\sigma^2(Y_t)\,dt$.
Since $\sigma^2=e^{2Y+2\theta}$ with $Y$ Gaussian, $\E[e^{Q}]=+\infty$. Hence $(S_T-K)^+$ has infinite
variance. A direct estimator still converges by the law of large numbers, but its confidence interval is meaningless.

\paragraph{Conditional (mixing) estimator.}
With $\rho=0$, $\log S_T$ given the path of $Y$ is Gaussian with variance $Q$, so
\begin{gather*}
\E\big[(S_T-K)^+\,\big|\,Y\big] = C_{\rm BS}\big(\sqrt{Q}\big),\\
C_{\rm BS}(s)=F\,\Phi(d_1)-K\,\Phi(d_1-s),\qquad
F=S_0e^{rT},\quad d_1=\frac{\ln(F/K)+s^2/2}{s}.
\end{gather*}
We average $C_{\rm BS}(\sqrt Q)$ over simulated paths of $Y$. Since $0\le C_{\rm BS}\le F$, the estimator has finite
variance and a valid confidence interval. Only $Y$ needs to be simulated.

\paragraph{Discretisation.}
$Y$ is advanced with the exact transition
$Y_{k+1}=m_1+a(Y_k-m_1)+\nu\sqrt{1-a^2}\,Z_k$, $a=e^{-\Delta t/\eps}$, $Z_k\sim\mathcal N(0,1)$. This step introduces no bias,
whatever the ratio $\Delta t/\eps$. $Q$ is computed with the trapezoidal rule and
$n=\max(250,\lceil 20T/\eps\rceil)$ steps, so that $\Delta t\le\eps/20$ and the quadrature resolves the fast scale.
We use $N=4\times10^4$ antithetic pairs ($Z\to-Z$). All strikes $K\in\{50,55,\dots,200\}$ are priced
on the same paths.
Implied volatilities are obtained by bisection on $C_{\rm BS}$.
The C++ code runs multithreaded, and the full sweep of $14$ values of $\eps$ takes about $5$\,s.

\paragraph{Remark on the Euler scheme.}
An Euler step for $Y$ has invariant variance $\nu^2/(1-\Delta t/(2\eps))$ instead of $\nu^2$, which inflates
$\E[\sigma^2]$ and hence the price. At $\eps=0.05$, $K=110$, Euler gives $55.94$ ($n=1000$) and $55.56$ ($n=4000$).
The exact transition gives $55.41\pm0.10$.

\section{Results}
\begin{table}[h]
\centering
\begin{tabular}{rrcccc}
\toprule
$\eps$ & $n$ & $\E[(S_T-110)^+]$ & 95\% CI & $e^{-rT}\E[\cdot]$ & implied vol\\
\midrule
1     & 250   & 35.027 & $\pm0.071$ & 32.334 & 0.848\\
0.5   & 250   & 39.390 & $\pm0.083$ & 36.361 & 0.959\\
0.2   & 250   & 46.087 & $\pm0.098$ & 42.544 & 1.136\\
0.1   & 250   & 51.228 & $\pm0.103$ & 47.290 & 1.277\\
0.05  & 400   & 55.406 & $\pm0.097$ & 51.146 & 1.397\\
0.02  & 1000  & 59.240 & $\pm0.082$ & 54.685 & 1.512\\
0.01  & 2000  & 61.077 & $\pm0.067$ & 56.381 & 1.569\\
0.005 & 4000  & 62.192 & $\pm0.053$ & 57.410 & 1.604\\
0.002 & 10000 & 62.952 & $\pm0.037$ & 58.112 & 1.628\\
0.001 & 20000 & 63.254 & $\pm0.028$ & 58.391 & 1.638\\
\midrule
0 (BS, $\bar\sigma$) & -- & 63.602 & -- & 58.713 & 1.649\\
\bottomrule
\end{tabular}
\caption{Strike $K=110$, $\rho=0$, $N=4\times10^4$ antithetic pairs.}
\label{tab:K110}
\end{table}

\begin{figure}[h]
\centering
\includegraphics[width=\textwidth]{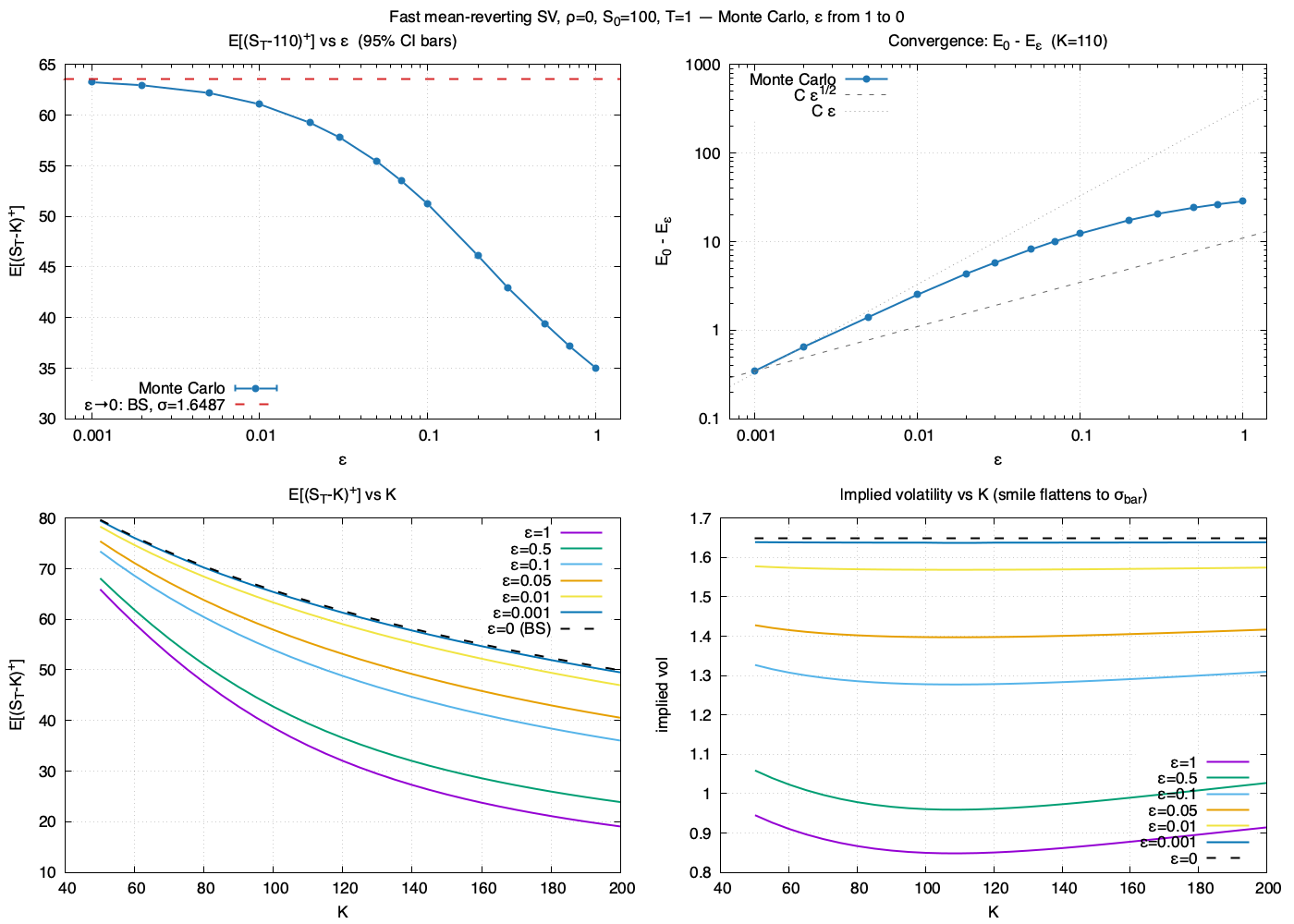}
\caption{Top left: $\E[(S_T-110)^+]$ versus $\eps$, with the Black--Scholes limit (dashed).
Top right: the gap $E_0-E_\eps$ in log--log scale, with reference slopes $\eps^{1/2}$ and $\eps$.
Bottom: price and implied volatility versus $K$ for $\eps\in\{1,0.5,0.1,0.05,0.01,0.001\}$ and $\eps=0$.}
\label{fig:eps}
\end{figure}

Table~\ref{tab:K110} and Figure~\ref{fig:eps} show three things.
\begin{enumerate}
\item The price increases monotonically with $1/\eps$ and converges to the Black--Scholes value $63.60$ with volatility $\bar\sigma$.
For large $\eps$, $Y$ (started at $m_1$) has not yet spread over its invariant law during $[0,T]$. Its variance is
$\nu^2(1-e^{-2t/\eps})$, so the realised variance, and hence the price, is much smaller.
\item The gap $E_0-E_\eps$ behaves like $O(\eps)$, not $O(\sqrt\eps)$. This is consistent with the
theory of Fouque, Papanicolaou and Sircar \cite{FPS}, \cite{FC}, (see also \cite{GM}), where the $\sqrt\eps$ correction is proportional to $\rho$ and vanishes here.
Both the transient and the fluctuations of $Q$ are $O(\eps)$. For example,
\[
T\bar\sigma^2-\E[Q] \;\approx\; \bar\sigma^2\,\frac{\eps}{2}\int_0^{2\nu^2}\frac{1-e^{-x}}{x}\,dx
\;\approx\; 0.66\,\bar\sigma^2\,\eps ,\qquad
\operatorname{Var}(Q)=O(\eps).
\]
\item For $\eps=O(1)$ the implied volatility is a smile, with its minimum near the money, as expected for a
volatility independent of $W^0$. The smile flattens and rises to the constant $\bar\sigma$ as $\eps\to0$.
\end{enumerate}

\section{A martingale control variate}\label{sec:cv}
We now simulate $S^\eps$ itself (log-exact step with $\sigma$ frozen on each step, $Y$ exact) and compare
\[
I^{st}=\frac1N\sum_{n=1}^N X^n,\qquad
I^{cv}=\frac1N\sum_{n=1}^N \big(X^n-M_T^n\big),\qquad X=e^{-rT}(S^\eps_T-K)^+ ,
\]
where $M_T$ is a stochastic integral with zero mean, computed on the same path.
Let $C(t,x)$ be the Black--Scholes price with volatility $\bar\sigma$, i.e.\ the solution of
$\partial_tC+rx\partial_xC+\frac12\bar\sigma^2x^2\partial_{xx}C-rC=0$, $C(T,x)=(x-K)^+$.
Its delta is $\partial_xC=\Phi(d_1)$.
We consider two martingales:
\[
M_T^{\bar\sigma}=\bar\sigma\int_0^Te^{-rt}\,\partial_xC(t,S_t)\,S_t\,dW_t,\qquad
M_T^{Y}=\int_0^Te^{-rt}\,\partial_xC(t,S_t)\,\sigma(Y_t)\,S_t\,dW_t ,
\]
where $W$ is the Brownian motion driving $S^\eps$. Both are discretised at left points, so
$\E[M_T]=0$ exactly and $I^{cv}$ is unbiased. Since $\operatorname{Var}(I)=\operatorname{Var}(X)/N$, we compare
variances per path.

\paragraph{Why $\sigma(Y_t)$ and not $\bar\sigma$.}
Itô's formula and the Black--Scholes equation give
\[
e^{-rT}(S_T-K)^+ = C(0,S_0) + M_T^{Y} + R_\eps,\qquad
R_\eps=\frac12\int_0^Te^{-rt}\big(\sigma^2(Y_t)-\bar\sigma^2\big)S_t^2\,\partial_{xx}C(t,S_t)\,dt .
\]
Hence $X-M_T^Y=C(0,S_0)+R_\eps$. Because $\langle\sigma^2\rangle=\bar\sigma^2$, the integrand of $R_\eps$ oscillates on the
time scale $\eps$ around zero mean, so $\operatorname{Var}(R_\eps)=O(\eps)$. Moreover $x^2\partial_{xx}C\le c/\sqrt{T-t}$, so this variance is finite,
whereas $\operatorname{Var}(X)=+\infty$ (Section~2).
With $\bar\sigma$ instead,
\[
X-M_T^{\bar\sigma} = C(0,S_0)+R_\eps+\int_0^Te^{-rt}\,\partial_xC\,\big(\sigma(Y_t)-\bar\sigma\big)S_t\,dW_t .
\]
The last term does not vanish as $\eps\to0$: its variance tends to
$\langle(\sigma-\bar\sigma)^2\rangle\,A$ with $A=\E\int_0^Te^{-2rt}(\partial_xC)^2S_t^2dt$.
Heuristically, $\operatorname{Var}(X)\approx\langle\sigma^2\rangle A$, and the best multiple
$\beta^*M^{\bar\sigma}$ has $\beta^*=\langle\sigma\rangle/\bar\sigma=e^{-1/2}\approx0.61$. It reduces the
variance by at most $\big(1-\langle\sigma\rangle^2/\langle\sigma^2\rangle\big)^{-1}=(1-e^{-1})^{-1}\approx1.6$.
Numerically we find $\beta^*\approx0.52$--$0.55$ and a factor $\approx1.4$ for $\eps\le10^{-2}$. With $\beta=1$ the factor is $\approx1$.

\begin{figure}[h]
\centering
\includegraphics[width=\textwidth]{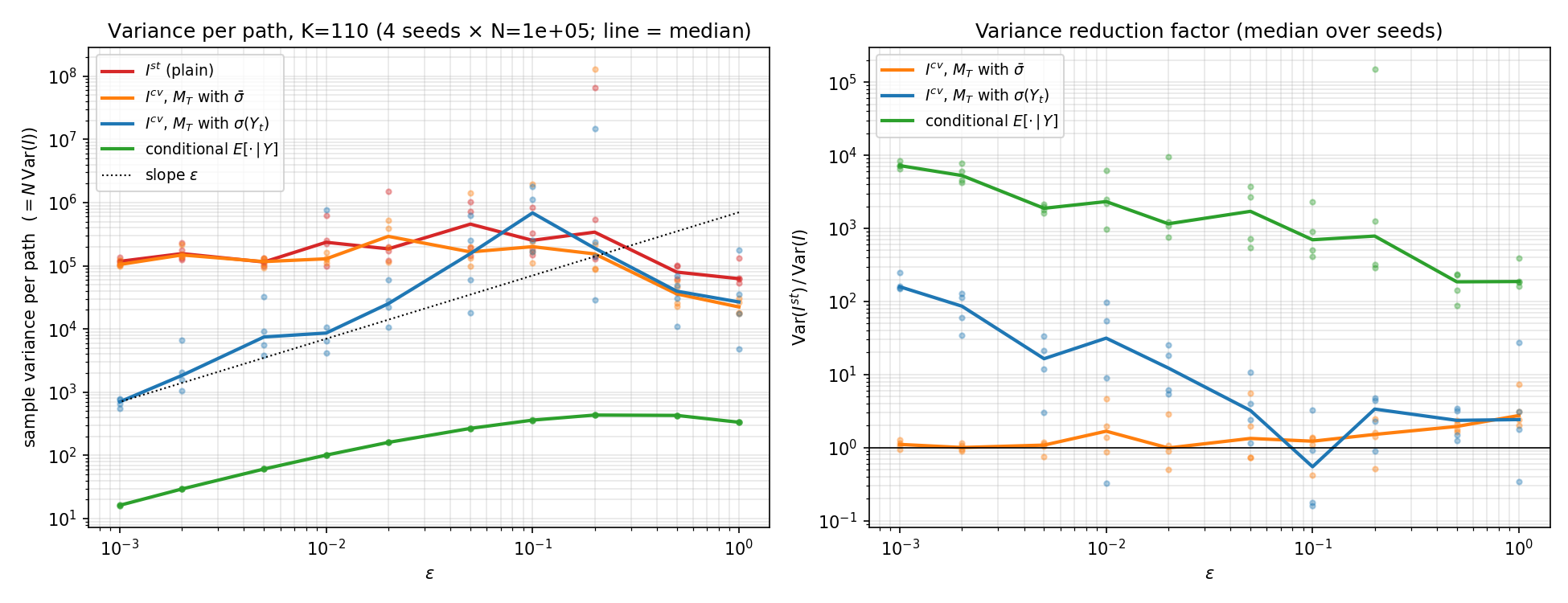}
\caption{Left: variance per path versus $\eps$ (dots: $4$ runs, line: median).
Right: reduction factor $\operatorname{Var}(I^{st})/\operatorname{Var}(I)$.}
\label{fig:cv}
\end{figure}

\paragraph{Results.}
We use $K=110$ and $4$ independent runs of $N=10^5$ paths for each $\eps$.
$\operatorname{Var}(X)$ is infinite, so its sample estimate is dominated by rare paths and changes from run to run.
We therefore report medians over the $4$ runs (Table~\ref{tab:cv}, Figure~\ref{fig:cv}).

\begin{table}[h]
\centering\small
\begin{tabular}{r|cccc|ccc|cc}
\toprule
& \multicolumn{4}{c|}{mean ($4\times10^5$ paths)} & \multicolumn{3}{c|}{variance per path (median)} &
\multicolumn{2}{c}{$\operatorname{Var}(I^{st})/\operatorname{Var}(I^{cv})$}\\
$\eps$ & $I^{st}$ & $I^{cv}_{\bar\sigma}$ & $I^{cv}_{Y}$ & cond. & $X$ & $X-M^{\bar\sigma}$ & $X-M^{Y}$ & $\bar\sigma$ & $\sigma(Y)$\\
\midrule
1 & 30.97 & 31.16 & 31.65 & 32.29 & 6.24e4 & 2.23e4 & 2.67e4 & 2.75 & 2.4\\
0.5 & 33.77 & 33.97 & 35.91 & 36.31 & 7.94e4 & 3.58e4 & 3.96e4 & 1.95 & 2.4\\
0.2 & 46.17 & 49.30 & 45.11 & 42.56 & 3.40e5 & 1.54e5 & 1.90e5 & 1.51 & 3.4\\
0.1 & 45.61 & 46.59 & 48.34 & 47.22 & 2.52e5 & 2.00e5 & 6.85e5 & 1.22 & 0.5\\
0.05 & 50.65 & 47.82 & 50.04 & 51.18 & 4.56e5 & 1.66e5 & 1.56e5 & 1.33 & 3.2\\
0.02 & 54.53 & 53.88 & 54.34 & 54.70 & 1.85e5 & 2.92e5 & 2.51e4 & 0.99 & 12.3\\
0.01 & 55.74 & 55.54 & 57.00 & 56.32 & 2.35e5 & 1.29e5 & 8.62e3 & 1.68 & 31.5\\
0.005 & 55.75 & 57.00 & 57.24 & 57.36 & 1.15e5 & 1.16e5 & 7.48e3 & 1.08 & 16.4\\
0.002 & 58.12 & 57.67 & 58.07 & 58.10 & 1.56e5 & 1.48e5 & 1.84e3 & 1.00 & 86.1\\
0.001 & 58.03 & 59.04 & 58.42 & 58.38 & 1.17e5 & 1.05e5 & 7.04e2 & 1.11 & 158.9\\
\bottomrule
\end{tabular}
\caption{Martingale control variate, discounted values, $K=110$. ``cond.'' is the conditional estimator of Section~2
(variance per path $16$ at $\eps=10^{-3}$). The Black--Scholes limit is $58.71$.}
\label{tab:cv}
\end{table}

The numerical results show the following.
\begin{itemize}
\item With $M_T^{Y}$, $\operatorname{Var}(I^{cv})\ll\operatorname{Var}(I^{st})$ for small $\eps$. The variance per path decreases like $\eps$
($7\times10^2$ at $\eps=10^{-3}$ against $1.2\times10^5$), and the reduction factor reaches about $160$ at $\eps=10^{-3}$.
The optimal coefficient is then $\beta^*=1.00$, as predicted by the identity $X-M_T^Y=C(0,S_0)+R_\eps$.
For $\eps\gtrsim0.05$ the gain is only $0.5$--$3$, because $R_\eps$ is then not small. At these values of $\eps$ the sample variances are also dominated by the heavy tail.
\item With $M_T^{\bar\sigma}$, as in the formula above, the ratio stays $\approx1$ for all $\eps\le0.02$.
The control variate does not work in this form. The martingale must use the volatility that actually drives $S^\eps$.
In the terminology of Fouque and Han \cite{FH}, $M_T^Y$ is their martingale $M_0(P_{BS})$ built on the leading-order
price $V_0=C(t,x)$ of the expansion of \cite{FPS}; higher-order terms of the expansion are tested in Section~\ref{sec:fh}.
\item For $\rho=0$ the conditional estimator is better still: its variance per path is $16$ at $\eps=10^{-3}$, about $7000$ times smaller than for $I^{st}$.
The martingale control variate remains the natural tool when $\rho\neq0$, where no closed-form conditional price exists.
\end{itemize}

\section{Calibration to S\&P~500 options}\label{sec:calib}
We now calibrate the model, with $\rho\neq0$, to market implied volatilities, and compare the Monte Carlo
estimator with and without the martingale control variate of Section~\ref{sec:cv} at the same number of paths.

\begin{figure}[h]
\centering
\includegraphics[width=\textwidth]{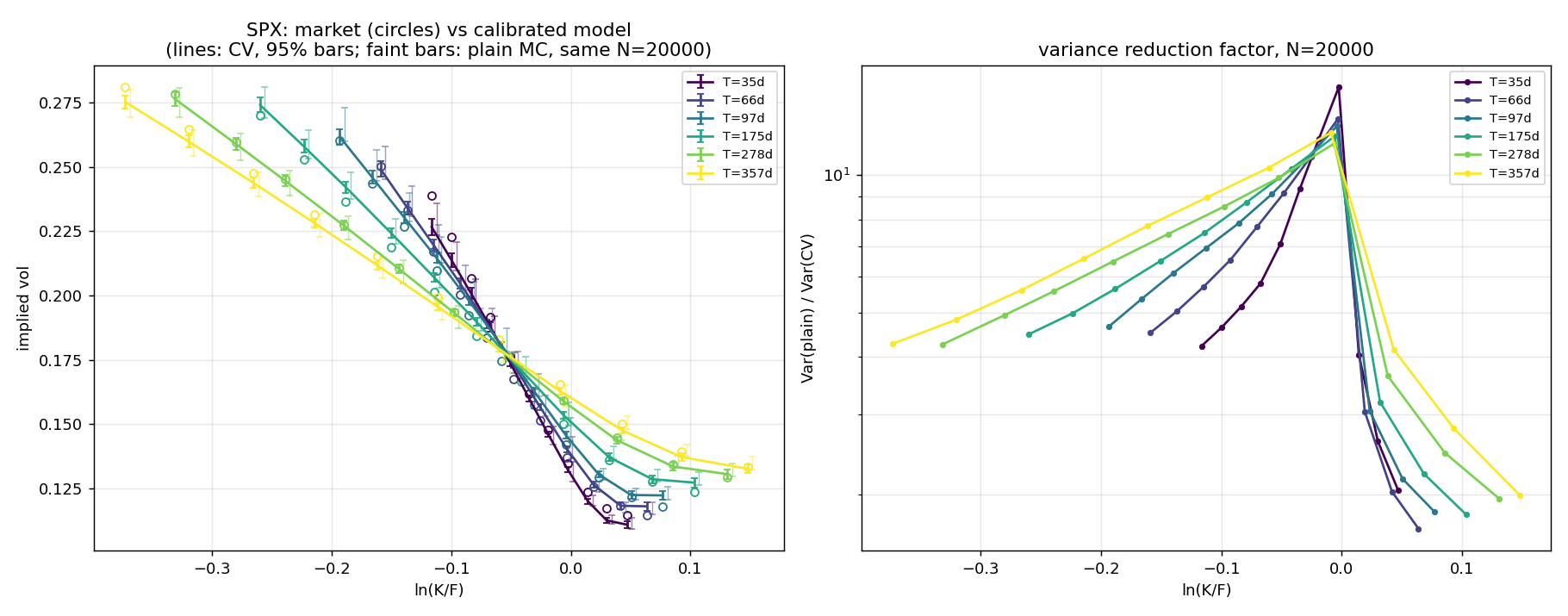}
\caption{Left: market implied volatilities (circles) and the calibrated model at the reference parameters,
with $95\%$ Monte Carlo intervals for $N=2\times10^4$ paths: with the control variate (solid bars) and without it
(faint bars, same paths). Right: variance reduction factor $\operatorname{Var}(X^{st})/\operatorname{Var}(X^{cv})$ per quote.}
\label{fig:calib}
\end{figure}

\paragraph{Market data.}
We use the S\&P~500 index options (ticker \texttt{\^{}SPX}, European) downloaded from Yahoo Finance on
25 September 2026 (spot $7704$). For each expiry the forward $F_T$ and the discount factor $D_T$ are obtained
from put--call parity, $C-P=D_T(F_T-K)$, by a linear regression over the $8$ strikes nearest to the money,
so no assumption on rates or dividends is needed. Implied volatilities are computed from bid--ask mid prices
of out-of-the-money options (puts for $K<F_T$, calls for $K\ge F_T$); quotes without a bid or with a spread
larger than half the mid price are discarded. We retain $6$ maturities ($35$, $66$, $97$, $175$, $278$ and $357$ days)
and, for each of them, the $11$ strikes nearest to $F_T e^{0.15\,z\sqrt T}$, $z\in\{-2.5,-2.15,\dots,1\}$,
i.e.\ $66$ quotes with implied volatilities from $0.11$ to $0.28$ (Figure~\ref{fig:calib}).

\paragraph{Model and parametrisation.}
Since $\sigma=e^{Y+\theta}$ with $\E[Y]=m_1$, the parameters $m_1$ and $\theta$ are confounded; we use the
log-volatility $X=Y+\theta$ and its mean $c=m_1+\theta$. Under the $T$-forward measure, with deterministic rates,
\[
dX_t=\frac1\eps(c-X_t)\,dt+\nu\sqrt{\frac2\eps}\,dW^1_t,\qquad
dS_t=S_t\,e^{X_t}\,dW^0_t,\quad S_0=1,\qquad d\langle W^0,W^1\rangle_t=\rho\,dt,
\]
and the call of maturity $T$ and strike $K$ is $D_TF_T\,\E[(S_T-k)^+]$ with $k=K/F_T$.
All maturities are priced on the same paths. The free parameters are $(\eps,c,\nu,\rho,X_0)$.

\paragraph{Estimators.}
For each quote, on each path, the plain estimator is the out-of-the-money payoff plus put--call parity,
\[
X^{st}=F_T\,(S_T-k)^+\ \ (k\ge1),\qquad X^{st}=F_T\,(k-S_T)^++F_T-K\ \ (k<1),
\]
which is unbiased because $\E[S_T]=1$ exactly in the scheme below; using the in-the-money call payoff
instead would inflate the variance of the plain estimator for low strikes. The control-variate estimator is
$X^{cv}=F_T(S_T-k)^+-M_T$ with the discrete hedging martingale
\[
M_T=F_T\sum_{b}\Delta_b\,\big(S_{t_{b+1}}-S_{t_b}\big),\qquad
\Delta_b=\Phi\Big(\frac{-\ln k}{s_b}+\frac{s_b}2\Big),\qquad
s_b^2=\int_{t_b}^T\E[\sigma_s^2]\,ds ,
\]
where $t_b$ are the rebalancing dates and $\E[\sigma_s^2]=\exp\big(2c+2(X_0-c)e^{-s/\eps}+2\nu^2(1-e^{-2s/\eps})\big)$
is known in closed form. $\Delta_b$ is the Black--Scholes delta with the deterministic variance term structure of the model,
which generalises $\bar\sigma$ of Section~\ref{sec:cv} to $X_0\neq c$. Since $S$ is an exact discrete martingale and
$\Delta_b$ is $\mathcal F_{t_b}$-measurable, $\E[M_T]=0$ for any choice of rebalancing dates, so $X^{cv}$ is unbiased.
When the delta is rebalanced at every step, $M_T\to M_T^Y$ as $\Delta t\to0$, and the identity of Section~\ref{sec:cv},
$X-M_T^Y=C(0,S_0)+R_\eps$, still holds for $\rho\neq0$ because $C$ depends only on $(t,S)$; the remainder is now
$R_\eps=\frac12\int_0^T(\sigma_t^2-\E[\sigma_t^2])S_t^2\partial_{xx}C\,dt$.
Note also that $X^{cv}$ is the same whether the call or the put is used, since both hedged payoffs differ by a constant.
Using the constant $\bar\sigma=e^{c+\nu^2}$ in $\Delta_b$ gives almost the same median variance reduction,
but a smaller one for short-dated out-of-the-money calls ($1.2$ instead of $1.6$).

\paragraph{Discretisation.}
$X$ is advanced with its exact transition $X_{k+1}=c+a(X_k-c)+\nu\sqrt{1-a^2}\,Z_k$, $a=e^{-\Delta t/\eps}$.
The Brownian increment $\Delta W^1_k$ is drawn jointly with the innovation $\xi_k=\nu^{-1}\sqrt{\eps/2}\,\sqrt{1-a^2}\,Z_k$
of the Ornstein--Uhlenbeck process, using $\operatorname{Var}\xi_k=\eps(1-a^2)/2$ and
$\operatorname{Cov}(\xi_k,\Delta W^1_k)=\eps(1-a)$, so that the correlation $\rho$ is exact for any ratio $\Delta t/\eps$.
Then $\Delta W^0_k=\rho\,\Delta W^1_k+\sqrt{1-\rho^2}\,\Delta W^\perp_k$ and
$\ln S_{k+1}=\ln S_k-\frac12\sigma_k^2\Delta t+\sigma_k\Delta W^0_k$ with $\sigma_k=e^{X_k}$.
We use $n=\lceil T_{\max}\max(250,4/\eps)\rceil$ steps.

\paragraph{Calibration.}
We minimise $\sum_q(\sigma^{\rm model}_q-\sigma^{\rm mkt}_q)^2$ over the $66$ quotes, where $\sigma^{\rm model}_q$ is the
Black implied volatility of the Monte Carlo price, with the trust-region reflective least-squares method of SciPy on the
unconstrained variables $(\ln\eps,c,\ln\nu,\operatorname{artanh}\rho,X_0)$, finite differences with relative step $10^{-2}$,
and common random numbers (the seed is fixed during a calibration), which makes the objective smooth in the parameters.
Starting from $(\eps,c,\nu,\rho,X_0)=(0.02,-2,0.5,-0.7,-2.1)$, each calibration takes $13$--$18$ pricings of the
$66$ quotes with $N=2\times10^4$ paths. We run it $5$ times (seeds $1$--$5$) with each estimator. A reference calibration
uses the control variate with $N=2\times10^5$. Each calibrated parameter set is then scored by its root-mean-square
implied-volatility error computed with an independent pricing at $N=2\times10^5$ (``true RMSE'').

\begin{table}[h]
\centering\small
\begin{tabular}{lcccccc}
\toprule
& $\eps$ & $c$ & $\nu$ & $\rho$ & $X_0$ & true RMSE (bp)\\
\midrule
reference ($N=2\times10^5$, CV) & $0.1956$ & $-2.4437$ & $0.9474$ & $-0.6515$ & $-2.1873$ & $35.8$\\
\midrule
with CV, mean of 5 & $0.1953$ & $-2.4491$ & $0.9473$ & $-0.6525$ & $-2.1848$ & $36.9$ \ ($36.0$--$38.6$)\\
with CV, RMS error & $0.0100$ & $0.0162$ & $0.0080$ & $0.0057$ & $0.0085$ & \\
without CV, mean of 5 & $0.2135$ & $-2.4237$ & $0.9436$ & $-0.6587$ & $-2.1803$ & $38.9$ \ ($36.6$--$41.2$)\\
without CV, RMS error & $0.0234$ & $0.0218$ & $0.0069$ & $0.0093$ & $0.0107$ & \\
\bottomrule
\end{tabular}
\caption{Calibrated parameters, $N=2\times10^4$ paths per pricing, $5$ seeds per estimator.
``RMS error'' is the root-mean-square deviation from the reference parameters.}
\label{tab:calib}
\end{table}

\paragraph{Calibrated model.}
The fit (Figure~\ref{fig:calib}, Table~\ref{tab:calib}) has a root-mean-square error of $36$ basis points of implied
volatility. The largest misfit, about $1.2$ vol points, is on the deep out-of-the-money puts of the shortest maturity,
whose skew is steeper than the model can produce. The correlation $\rho\approx-0.65$ generates the skew and the initial
volatility $e^{X_0}\approx0.112$, below the long-run level $\bar\sigma=e^{c+\nu^2}\approx0.213$, generates the
increasing term structure. The mean-reversion time is $\eps\approx0.2$ years, i.e.\ about $70$ days: the market
surface does not select the fast regime $\eps\ll1$ for this one-factor model.

\paragraph{Precision with and without the control variate.}
At the reference parameters and with the same $N=2\times10^4$ paths, the control variate reduces the variance by a
factor $1.7$--$15.6$ (median $5.6$). The median standard error in implied volatility goes from $22.3$ to $8.0$ basis points.
The gain is largest at the money ($12$--$16$), where the delta hedge removes most of the payoff variance; it is $4$--$5$
for the out-of-the-money puts and only $1.7$--$2$ for the out-of-the-money calls, whose variance comes mostly from
volatility paths and is not hedged by a delta in $S$. At $\eps\approx0.2$ this is consistent with Section~\ref{sec:cv}, where
the gain was a few units and grew like $1/\eps$ only for $\eps\lesssim0.02$.

The smaller Monte Carlo noise gives better calibrated parameters (Table~\ref{tab:calib}). With the control variate the
five calibrations have a true RMSE of $36.0$--$38.6$\,bp, against $36.6$--$41.2$\,bp without it (reference: $35.8$\,bp),
and the RMS deviation from the reference parameters is smaller for $\eps$ ($0.010$ against $0.023$), $c$, $\rho$ and $X_0$;
it is similar for $\nu$. Without the control variate, the mean of $\eps$ is also $0.018$ above the reference, comparable to its scatter ($0.014$),
so five runs do not establish a bias.

\paragraph{Cost.}
The price of the control variate is the evaluation of $\Delta_b$ for every quote at every rebalancing date.
Table~\ref{tab:hedge} gives the variance reduction and the CPU time relative to the plain estimator, for $N=2\times10^5$
and several numbers of rebalancing dates per year. With a rebalancing at every step ($250$ per year) the control variate
costs $10.6$ times more than plain Monte Carlo, so it is less efficient at equal CPU time (efficiency $0.52$).
With a coarser rebalancing the martingale remains exactly centred, the reduction decreases slowly and the cost
decreases quickly: with $25$ dates per year the efficiency is $1.7$.

\begin{table}[h]
\centering
\begin{tabular}{rcccc}
\toprule
dates / year & \multicolumn{2}{c}{variance ratio: median (range)} & time CV / time plain & efficiency\\
\midrule
$250$ (every step) & $5.50$ & ($1.59$--$15.1$) & $10.6$ & $0.52$\\
$100$ & $4.69$ & ($1.44$--$12.1$) & $4.57$ & $1.03$\\
$50$  & $4.11$ & ($1.27$--$11.1$) & $3.00$ & $1.37$\\
$25$  & $3.29$ & ($1.11$--$9.59$) & $1.95$ & $1.69$\\
\bottomrule
\end{tabular}
\caption{Martingale control variate at the reference parameters, $N=2\times10^5$ paths, $66$ quotes.
Efficiency $=$ median variance ratio $/$ time ratio; values above $1$ favour the control variate at equal CPU time.}
\label{tab:hedge}
\end{table}

\begin{figure}[h]
\centering
\includegraphics[width=\textwidth]{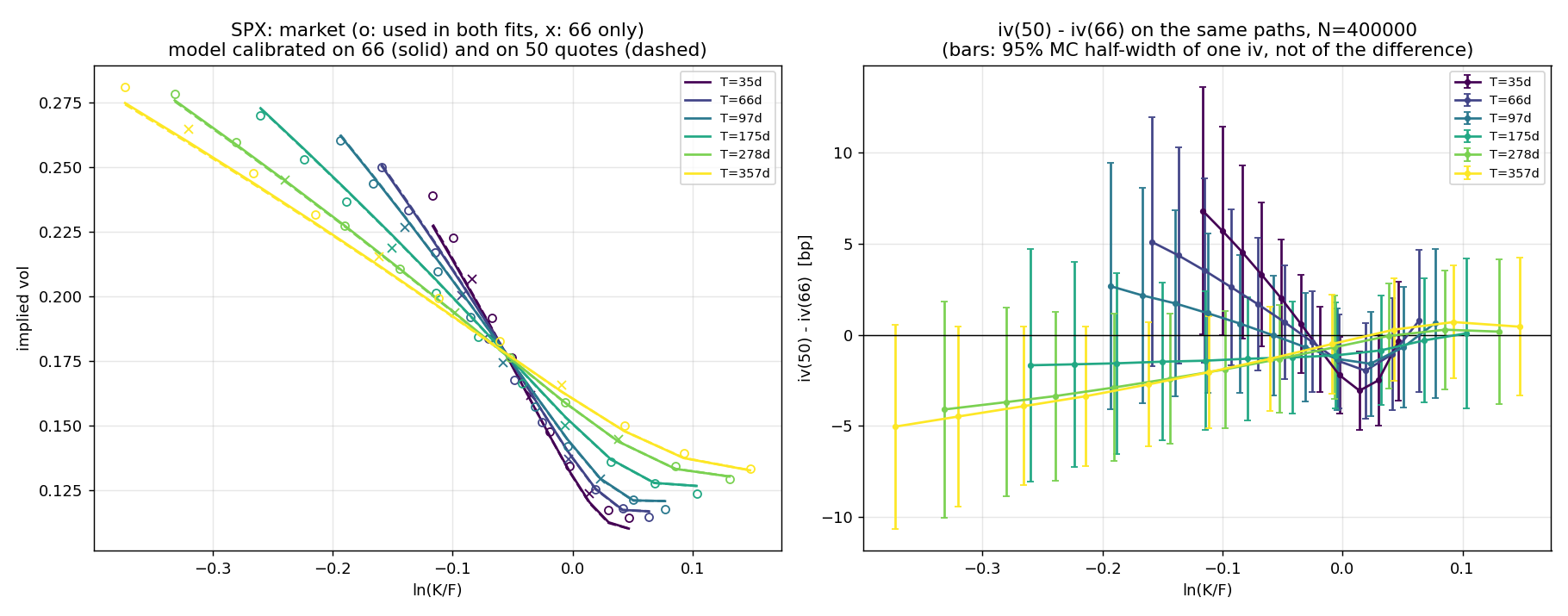}
\caption{Left: market implied volatilities (circles: used in both calibrations; crosses: used only with $66$ quotes)
and the model calibrated on $66$ (solid) and on $50$ quotes (dashed). Right: difference of the two model implied
volatilities in basis points, computed on the same $4\times10^5$ paths; the bars are the $95\%$ Monte Carlo half-width
of one implied volatility, an upper bound for the noise of the difference.}
\label{fig:stab}
\end{figure}

\paragraph{Stability with respect to the data.}
To test how much the result depends on the quotes, we calibrate again on a subset of $50$ of the $66$ quotes.
In each maturity we remove $2$ or $3$ interior strikes and keep the two extreme ones, so that the smaller calibration
never extrapolates. Both calibrations use the control variate, $N=2\times10^5$ paths, the same seed and the same
starting point, so they differ only by the data. Both parameter sets are then priced on all $66$ quotes with the same
$4\times10^5$ paths (Table~\ref{tab:stab}, Figure~\ref{fig:stab}).

\begin{table}[h]
\centering\small
\begin{tabular}{lccccc|ccc}
\toprule
& & & & & & \multicolumn{3}{c}{RMSE vs market (bp)}\\
quotes used & $\eps$ & $c$ & $\nu$ & $\rho$ & $X_0$ & 50 common & 16 others & all 66\\
\midrule
50 & $0.1902$ & $-2.4491$ & $0.9486$ & $-0.6505$ & $-2.1922$ & $37.3$ & $35.5$ & $36.8$\\
66 & $0.1961$ & $-2.4436$ & $0.9474$ & $-0.6516$ & $-2.1868$ & $37.4$ & $35.1$ & $36.8$\\
\bottomrule
\end{tabular}
\caption{Calibration on $50$ and on $66$ quotes ($N=2\times10^5$, control variate). For the calibration on $50$ quotes,
the $16$ other quotes are out of sample.}
\label{tab:stab}
\end{table}

The calibration is stable. The parameter $\eps$ changes by $3\%$ and the others by less than $0.3\%$, much less than the
scatter due to Monte Carlo noise at $N=2\times10^4$ in Table~\ref{tab:calib}. The two implied-volatility surfaces differ
by $2.5$\,bp in root-mean-square value and by at most $6.8$\,bp; this is more than ten times smaller than the misfit of
the model ($37$\,bp). The largest differences are on the out-of-the-money puts of the shortest maturity, where the model
also fits the market least well; beyond $175$ days they are below $2$\,bp. Finally, the calibration on $50$ quotes
prices the $16$ quotes it has not seen as well as the calibration that uses them ($35.5$ against $35.1$\,bp), so the
$66$ quotes do not overfit the model: the remaining error is a model error, not a sampling effect of the quotes.

\subsection{Higher-order approximations in the martingale}\label{sec:fh}
The delta used so far is that of the leading-order term $V_0$ of the fast mean-reversion expansion of
Fouque, Papanicolaou and Sircar \cite{FPS}, as in Fouque and Han \cite{FH}. We test whether the next terms of the
expansion give a better control variate. With the operator $\mathcal L_0=(c-y)\partial_y+\nu^2\partial_{yy}$ of the
unscaled Ornstein--Uhlenbeck process, the pricing operator of the model of this section is
$\frac1\eps\mathcal L_0+\frac1{\sqrt\eps}\mathcal L_1+\mathcal L_2$ with $\mathcal L_1=\sqrt2\,\nu\rho\,\sigma(y)\,x\,\partial_{xy}$,
so the expansion is in powers of $\sqrt\eps$: $P^\eps=V_0+\sqrt\eps V_1+\eps V_2+\cdots$, with
\[
V_1=-(T-t)\,\mathcal C\,x\partial_x\big(x^2\partial_{xx}V_0\big),\qquad
\mathcal C=\frac{\nu\rho}{\sqrt2}\langle\sigma\phi'\rangle,\qquad
V_2=-\tfrac12\,\phi(y)\,x^2\partial_{xx}V_0+\bar V_2(t,x),
\]
where $V_0$ is the Black price with volatility $\bar\sigma$ and $\phi$ solves $\mathcal L_0\phi=\sigma^2-\bar\sigma^2$,
$\langle\phi\rangle=0$. For $\sigma=e^y$ and the invariant law $\mathcal N(c,\nu^2)$ both are explicit: an integration by parts
gives $\langle\sigma\phi'\rangle=-\big(\langle\sigma^3\rangle-\langle\sigma\rangle\bar\sigma^2\big)/\nu^2$, and the Hermite expansion
$e^{2y}=\bar\sigma^2\sum_{n\ge0}\frac{(2\nu)^n}{n!}\mathrm{He}_n(z)$, $z=(y-c)/\nu$, with $\mathcal L_0\mathrm{He}_n=-n\,\mathrm{He}_n$, gives
$\phi(y)=-\bar\sigma^2\sum_{n\ge1}\frac{(2\nu)^n}{n\,n!}\mathrm{He}_n(z)$. Note that $V_1=0$ when $\rho=0$, so $V_1$ cannot improve the
control variate of Section~\ref{sec:cv}.

We compare, on the same paths and with rebalancing at every step, the martingales built on
(A) $V_0$ with $\bar\sigma$; (B) $V_0$ with the variance term structure (the choice of this section);
(C) $V_0+\sqrt\eps V_1$; and (D) $V_0+\sqrt\eps V_1-\frac\eps2\phi(y)x^2\partial_{xx}V_0$. For (D) the approximation depends
on $y$, so the martingale also contains $\sum_k\partial_y\tilde P\,\big(X_{k+1}-\E[X_{k+1}\mid X_k]\big)$, the term $M_1$ of \cite{FH};
it has exactly zero mean with the exact Ornstein--Uhlenbeck step. The unknown $\bar V_2$ does not depend on $y$ and is omitted.

\begin{table}[h]
\centering\small
\begin{tabular}{lcccc}
\toprule
& \multicolumn{3}{c}{median variance ratio} & IV error (bp)\\
& A: $V_0$ & C: $+\sqrt\eps V_1$ & D: C $+\,\eps V_2^\perp$ & $V_0$ / $V_0+\sqrt\eps V_1$\\
\midrule
SPX ($\eps=0.196$, $\nu=0.95$) & $6.1$ & $0.96$ & $0.03$ & $563$ / $1633$\\
$\nu=0.3$, $\eps=0.2$   & $27$  & $17$  & $5.2$ & $126$ / $55$\\
$\nu=0.3$, $\eps=0.01$  & $96$  & $91$  & $159$ & $44$ / $9.1$\\
$\nu=0.3$, $\eps=0.002$ & $458$ & $444$ & $844$ & $20$ / $1.9$\\
$\nu=0.95$, $\eps=0.01$  & $12.6$ & $9.0$ & $6.7$ & $197$ / $122$\\
$\nu=0.95$, $\eps=0.002$ & $41.6$ & $33.9$ & $29.3$ & $94$ / $34$\\
\bottomrule
\end{tabular}
\caption{Median variance reduction $\operatorname{Var}(X^{st})/\operatorname{Var}(X^{cv})$ with the martingales built on successive
terms of the expansion. First row: calibrated SPX parameters, $66$ quotes, $N=2\times10^4$ (variant B gives $5.6$).
Other rows: $c$ and $\rho$ of SPX, $X_0=c$, $T\in\{0.5,1\}$, $k\in\{0.85,0.95,1,1.05,1.15\}$, $N=10^5$.
Last column: root-mean-square implied-volatility error of the approximate prices against Monte Carlo.}
\label{tab:fh}
\end{table}

\begin{figure}[h]
\centering
\includegraphics[width=\textwidth]{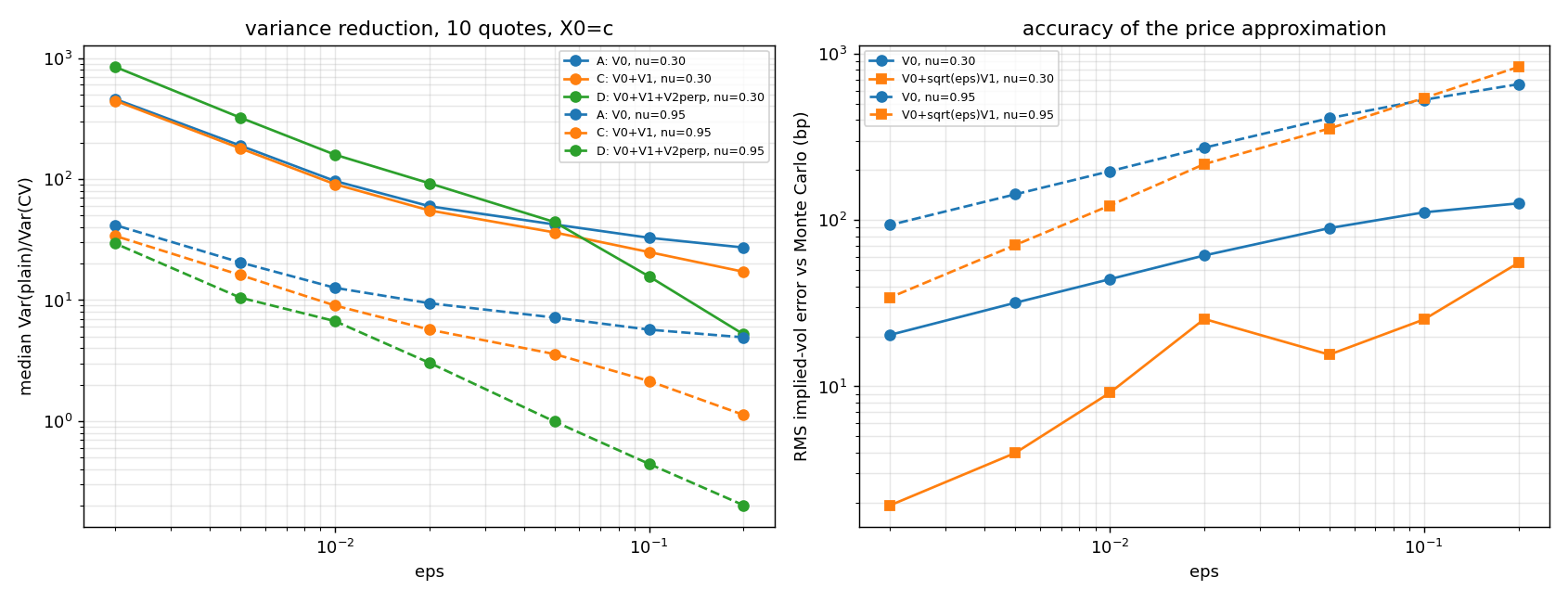}
\caption{Left: median variance reduction versus $\eps$ for the martingales built on $V_0$, $V_0+\sqrt\eps V_1$ and
$V_0+\sqrt\eps V_1+\eps V_2^\perp$, for $\nu=0.3$ (solid) and $\nu=0.95$ (dashed). Right: implied-volatility error of the
approximate prices $V_0$ and $V_0+\sqrt\eps V_1$.}
\label{fig:fh}
\end{figure}

The results (Table~\ref{tab:fh}, Figure~\ref{fig:fh}) are the following.
\begin{itemize}
\item The first-order price $V_0+\sqrt\eps V_1$ is much more accurate than $V_0$ when the expansion is valid: for $\nu=0.3$ its
error is $1.9$\,bp against $20$\,bp at $\eps=0.002$. This validates the formulas and the scaling.
\item Nevertheless $V_1$ never improves the control variate: C is slightly worse than A for small $\eps$ and much worse for large
$\eps$ or $\nu$. This agrees with the remark of Fouque and Han \cite{FH} that the next-order correction leaves the variance of the
same order: $V_1$ does not depend on $y$, whereas the residual variance comes from $\frac1\eps\int(\partial_yP^\eps)^2g^2\,dt$,
i.e.\ from the volatility factor, which a delta in $x$ cannot hedge.
\item The $y$-dependent part $V_2^\perp$, with the martingale $M_1$ in the volatility factor, addresses this term. It doubles the
variance reduction of $V_0$ when the expansion is valid ($844$ against $458$ for $\nu=0.3$, $\eps=0.002$), but it is harmful as soon as
$\eps\gtrsim0.05$ or $\nu$ is large, because $\phi(y)$ grows like $e^{2y}$ and $\eps\phi$ is then not small.
\item At the calibrated SPX parameters ($\eps\approx0.2$, $\nu\approx0.95$) the expansion is not valid: even the price $V_0+\sqrt\eps V_1$
is worse than $V_0$, C removes almost all the gain and D increases the variance about $30$ times. The best choice for market-calibrated
models remains $V_0$, preferably with the variance term structure.
\end{itemize}

\section{A basket of four indices}\label{sec:basket}
Monte Carlo is the method of choice when the dimension is large. We therefore price a call on a basket of four
US indices, each following the model of Section~\ref{sec:calib} with its own calibrated parameters. The state
has dimension $8$ (four prices and four volatility factors), which is out of reach of grid-based PDE solvers,
while the cost of the Monte Carlo simulation grows only linearly with the number of assets.

\paragraph{Marginals.}
The four indices are the S\&P~500 (SPX), the Dow Jones, the Nasdaq-100 (NDX) and the Russell~2000 (RUT).
The European options on the Dow (DJX) have too few quotes on Yahoo Finance, so we use the options on the DIA
exchange-traded fund; they are American, but we only use out-of-the-money mid prices, for which the early-exercise
premium is small. The chains of DIA, NDX and RUT contain stale quotes, which we remove before the selection of
Section~\ref{sec:calib}: a maturity is rejected when the put and call volatilities do not join at the forward
(gap larger than $1.5$ vol points, the sign of a wrong forward or of a stale side), quotes further than
$\max(1\%,3\,\mathrm{MAD})$ from a robust cubic fit of the smile are dropped, a maturity is rejected when the
remaining residuals have a standard deviation above $0.6$ vol points, and at least $8$ strikes are required per
maturity. We keep the maturities nearest to $35$, $84$, $175$ and $357$ days that pass these tests.
Each index is then calibrated separately with the one-dimensional pricer ($N=2\times10^5$, control variate);
SPX uses the reference calibration of Table~\ref{tab:calib}.

\begin{table}[h]
\centering\small
\begin{tabular}{lccccccc}
\toprule
index & $\eps$ & $c$ & $\nu$ & $\rho$ & $X_0$ & maturities (days) & RMSE (bp)\\
\midrule
SPX & $0.1956$ & $-2.4437$ & $0.9474$ & $-0.6515$ & $-2.1873$ & $35$--$357$ ($6$), $66$ quotes & $35.8$\\
DIA & $0.1182$ & $-2.3849$ & $0.8245$ & $-0.4236$ & $-2.1543$ & $35$, $84$, $22$ quotes & $14.9$\\
NDX & $0.1126$ & $-2.0464$ & $0.8265$ & $-0.4610$ & $-1.8155$ & $35$, $84$, $175$, $33$ quotes & $28.7$\\
RUT & $0.1537$ & $-2.0636$ & $0.7607$ & $-0.4812$ & $-1.7555$ & $35$, $84$, $126$, $265$, $44$ quotes & $42.3$\\
\bottomrule
\end{tabular}
\caption{Calibrated marginals of the basket (market data of 25 September 2026).}
\label{tab:marg}
\end{table}

\paragraph{Basket model.}
For each index $i=1,\dots,4$ let $S^i_t$ be its forward-normalised performance, a martingale with $S^i_0=1$ that follows
the model of Section~\ref{sec:calib} with the parameters of Table~\ref{tab:marg}. The Brownian motions of different indices
are independent; within one index, $d\langle W^{0,i},W^{1,i}\rangle=\rho_i\,dt$. The basket is
$B_t=\sum_iw_iS^i_t$ with $w_i=1/4$, so that $\E[B_T]=1$, and we price $\E[(B_T-k)^+]$ for $T=3$ and $6$ months and
$9$ strikes $k=e^{0.1z\sqrt T}$, $z\in[-2.5,1.5]$. The maturities are limited to $6$ months because the DIA marginal is
calibrated up to $84$ days only. Independence is a simplification chosen to keep the example transparent: US indices are
strongly correlated, so the basket volatility ($\approx0.10$) is about half of a realistic one. A correlation between
the indices only changes the Gaussian increments and the basket variance $s_b^2$ below.

\paragraph{Estimators.}
As in dimension one, the plain estimator is the out-of-the-money payoff plus parity,
$X^{st}=(B_T-k)^+$ for $k\ge1$ and $(k-B_T)^++1-k$ for $k<1$. The control variate hedges the basket with the Black delta
of a lognormal basket,
\begin{gather*}
X^{cv}=(B_T-k)^+-\sum_b\Delta_b\big(B_{t_{b+1}}-B_{t_b}\big),\\
\Delta_b=\Phi\Big(\frac{\ln(B_{t_b}/k)}{s_b}+\frac{s_b}2\Big),\qquad
s_b^2=\sum_i\Big(\frac{w_iS^i_{t_b}}{B_{t_b}}\Big)^2\int_{t_b}^T\E[(\sigma^i_s)^2]\,ds .
\end{gather*}
The martingale is the profit of holding $w_i\Delta_b$ units of each index; since every $S^i$ is an exact discrete martingale
and $\Delta_b$ is $\mathcal F_{t_b}$-measurable, it has zero mean and $X^{cv}$ is unbiased, whatever the quality of the
lognormal approximation. The simulation of each index is that of Section~\ref{sec:calib}. With weights $(1,0,0,0)$ the basket
code reproduces the one-dimensional pricer within one standard error, with identical standard errors.

\begin{figure}[h]
\centering
\includegraphics[width=\textwidth]{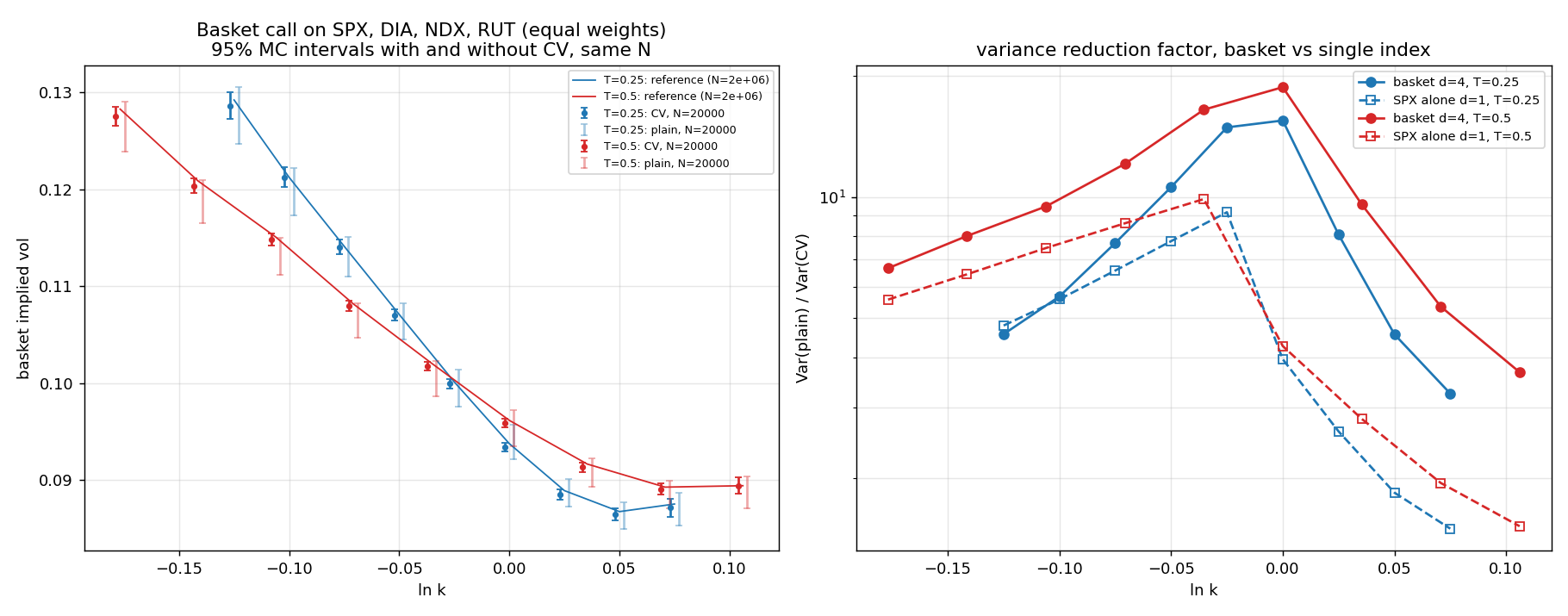}
\caption{Left: implied volatility of the basket call (reference: $N=2\times10^6$ with control variate) and $95\%$ Monte Carlo
intervals with $N=2\times10^4$ paths, with the control variate (solid) and without it (faint, same paths).
Right: variance reduction factor for the basket ($d=4$, solid) and for SPX alone ($d=1$, dashed) on the same quotes.}
\label{fig:basket}
\end{figure}

\begin{table}[h]
\centering
\begin{tabular}{crcccc}
\toprule
$d$ & dates / year & \multicolumn{2}{c}{variance ratio: median (range)} & time CV / time plain & efficiency\\
\midrule
$4$ & $250$ (every step) & $8.25$ & ($3.41$--$19.5$) & $2.25$ & $3.66$\\
$4$ & $100$ & $6.71$ & ($2.85$--$16.6$) & $1.30$ & $5.16$\\
$4$ & $50$  & $5.61$ & ($2.35$--$14.0$) & $1.25$ & $4.50$\\
$4$ & $25$  & $4.28$ & ($2.03$--$11.3$) & $0.98$ & $4.36$\\
\midrule
$1$ & $250$ (every step) & $5.71$ & ($1.62$--$10.6$) & $6.15$ & $0.93$\\
$1$ & $100$ & $4.92$ & ($1.47$--$9.59$) & $2.69$ & $1.83$\\
$1$ & $50$  & $4.46$ & ($1.36$--$8.95$) & $2.02$ & $2.20$\\
$1$ & $25$  & $3.62$ & ($1.21$--$7.64$) & $1.70$ & $2.14$\\
\bottomrule
\end{tabular}
\caption{Basket ($d=4$) and SPX alone ($d=1$) on the same $18$ quotes, $N=2\times10^5$ paths. Efficiency $=$ median variance
ratio $/$ time ratio. The run times are between $0.1$ and $1.1$\,s on $14$ threads and are indicative only.}
\label{tab:bhedge}
\end{table}

\paragraph{Results.}
With the same $N=2\times10^4$ paths, the control variate reduces the variance of the basket prices by a factor $3.3$--$18.8$
(median $8.0$), and the median standard error of the basket implied volatility goes from $9.4$ to $3.1$ basis points
(Figure~\ref{fig:basket}). On the same quotes, SPX alone gains only a factor $1.5$--$9.9$ (median $5.2$).
The two gains are equal for the lowest strikes at $3$ months; elsewhere the basket gains more, and much more at the money
($15$--$19$ instead of about $4$) and above it ($3.3$--$9.6$ instead of $1.5$--$2.8$).
A heuristic reason is that the remainder of the hedge, $R_\eps=\frac12\int_0^T\sum_iw_i^2(S^i_t)^2\big((\sigma^i_t)^2-\E[(\sigma^i_t)^2]\big)
\partial_{BB}C\,dt$, is a sum of independent fluctuations and averages out, whereas the variance of the payoff itself is
dominated by the common level of the basket, which the delta removes.

The cost is also more favourable (Table~\ref{tab:bhedge}). The simulation cost is proportional to $d$, but the hedge needs one
delta per quote and rebalancing date for the whole basket. With a rebalancing at every step, the control variate costs $2.25$
times the plain estimator for $d=4$ against $6.15$ for $d=1$, so it is $3.7$ times more efficient at equal CPU time for the
basket, while it only breaks even for a single index. With $100$ rebalancing dates per year the efficiency reaches $5.2$.
In dimension $4$ the martingale control variate is therefore worthwhile in its simplest form.

\section{Conclusion}
The conditional Monte Carlo estimator \cite{G} with exact simulation of the fast factor is accurate
($\pm0.03$--$0.1$ at 95\%) and cheap for all $\eps\in[10^{-3},1]$. It confirms the averaging
limit $\eps\to0$ and gives an $O(\eps)$ convergence rate when $\rho=0$.
The case $\rho\neq0$ requires $I_1=\int_0^T\sigma(Y_t)\,dW^1_t$ in the conditional forward. There the $\sqrt\eps$ correction
should reappear.
The control variate $X-M_T$, with $M_T$ built on the Black--Scholes delta at volatility $\bar\sigma$, gives
$\operatorname{Var}(I^{cv})=O(\eps)\ll\operatorname{Var}(I^{st})$, provided that $M_T$ is driven by the true volatility $\sigma(Y_t)$.
Calibrated to S\&P~500 options (Section~\ref{sec:calib}), the model with $\rho\approx-0.65$ fits $66$ quotes within
$36$ basis points of implied volatility, with $\eps\approx0.2$ years, outside the fast regime. There the martingale
control variate, written as a discrete delta hedge with the model's variance term structure, reduces the variance
$5.6$ times at the median for the same number of paths and improves the calibrated parameters. Its cost comes from the
deltas; rebalancing them about $25$ times a year keeps most of the gain at twice the cost of plain Monte Carlo.
The calibration is stable: removing $16$ of the $66$ quotes changes the calibrated implied volatilities by $2.5$\,bp
(root-mean-square), and the removed quotes are fitted as well as the others.
For a call on a basket of four indices (SPX, Dow, Nasdaq-100, Russell~2000; Section~\ref{sec:basket}) the same construction,
with the Black delta of the basket, reduces the variance $8$ times at the median and costs only about twice the plain
estimator, so that the method gains a factor $4$--$5$ at equal CPU time; this is where it is most useful.

\subsection*{Acknowledgement}
Most of the work has been done by Claude.Code, the AI engine of Anthropic ltd.Latex layout and write-up, Monte Carlo coding in C++ with gnuplot displays. Even the idea to replace $\bar\sigma$ by $\sigma(Y)$ for the variance estimator is Claude's idea.

The human authors asked the questions and monitored the answers.


\end{document}